# Existence and Uniqueness Theorems for Sequential Linear Conformable Fractional Differential Equations

Ahmet Gökdoğan[a], Emrah Ünal[b], Ercan Çelik[c]

[a] Department of Mathematical Engineering, Gümüşhane University, 29100 Gümüşhane, Turkey,

gokdogan@gumushane.edu.tr

[b] Department of Elementary Mathematics Education, Artvin Çoruh University, 08100 Artvin, Turkey

emrah.unal@artvin.edu.tr

[c] Department of Mathematics, Atatürk University, 25400 Erzurum, Turkey,

ecelik@atauni.edu.tr

**Abstract**

Recently, a new fractional derivative called the conformable fractional derivative is given on based basic limit definition derivative in [4]. Then, the fractional versions of chain rules, exponential functions, Gronwall's inequality, integration by parts, Taylor power series expansions is developed in [5]. In this paper, we present existence and uniqueness theorems for sequential linear conformable fractional differential equations.

**MSC:** 34AXX

**Key Words:** Sequential Linear Fractional Differential Equations, Conformable Fractional Derivative, Existence and Uniqueness Theorems

## 1. Introduction

Although the idea of fractional derivative has been suggested by L'Hospital at the last 17.Century, intensive studies about fractional derivative were carried out in the last and present centuries. Many researchers used an integral form for fractional derivative definition. The most popular definitions of fractional derivative are Riemann-Liouville and Caputo definitions. For Riemann-Liouville, Caputo and other definitions and the characteristics of these definitions, we refer to reader to [1-3].

(I) Riemann Liouville definition:
$$D_x^\alpha f(x) = \frac{1}{\Gamma(n-\alpha)} \left(\frac{d}{dx}\right)^n \int_0^x (x-t)^{n-\alpha-1} f(t) dt, \quad n-1 < \alpha \leq n$$

(II) Caputo definition:
$$D_x^\alpha f(x) = \frac{1}{\Gamma(n-\alpha)} \int_0^x (x-t)^{n-\alpha-1} \left(\frac{d}{dx}\right)^n f(t) dt, \quad n-1 < \alpha \leq n$$

Lately, Khalil at al. give a new definition of fractional derivative and fractional integral [4]. This new definition benefit from a limit form as in usual derivatives. This new theory is improved by Abdeljawad. For instance, Taylor power series representation and Laplace transform of few certain functions, fractional integration by parts formulas, chain rule and Gronwall inequality are provided by him [5]. Many studies about theory and application of the sequential linear fractional differential equations which based on Hadamard, Riemann-Liouville and Caputo derivatives were conducted [6-11].

In this study, we present existence and uniqueness theorems for sequential linear conformable fractional differential equations.

## 2. Conformable Fractional Derivative

Here, some basic definitions and properties of the conformable fractional calculus theory which can be found in [4–5].

**Definition 2.1.** $f:[0,\infty) \to \mathbb{R}$ let a function. Then for all $t > 0$, the conformable fractional derivative of $f$ of order $\alpha$ is defined as

$$T_\alpha(f)(t) = \lim_{\varepsilon \to 0} \frac{f(t+\varepsilon t^{1-\alpha}) - f(t)}{\varepsilon} \tag{1}$$

where $\alpha \in (0,1)$.

If $f$ is $\alpha$-differentiable in some $(0, a), a > 0$ and $\lim_{t\to 0^+} f^{(\alpha)}(t)$ exists, then define

$$f^{(\alpha)}(0) = \lim_{t \to 0^+} f^{(\alpha)}(t). \tag{2}$$

**Theorem 2.1.** Let $\alpha \in (0,1]$ and $f, g$ be $\alpha$-differentiable at a point $t > 0$. Then

(1) $T_\alpha(af + bg) = aT_\alpha(f) + bT_\alpha(g)$, for all $a, b \in \mathbb{R}$.

(2) $T_\alpha(t^p) = pt^{p-\alpha}$ for all $p \in \mathbb{R}$.

(3) $T_\alpha(\lambda) = 0$ for all constant functions $f(t) = \lambda$.

(4) $T_\alpha(fg) = T_\alpha(f)g + fT_\alpha(g)$.

(5) $T_\alpha(f/g) = \frac{T_\alpha(f)g - fT_\alpha(g)}{g^2}$.

(6) In addition, If $f$ is differentiable, then $T_\alpha(f(t)) = t^{1-\alpha}\frac{df}{dt}$.

Additionaly, conformable fractional derivatives of certain functions as follow:

(i) $\quad T_\alpha\left(\sin\frac{1}{\alpha}t^\alpha\right) = \cos\frac{1}{\alpha}t^\alpha$.

(ii) $T_\alpha \left(\cos \frac{1}{\alpha} t^\alpha\right) = -\sin \frac{1}{\alpha} t^\alpha.$

(iii) $T_\alpha \left(e^{\frac{1}{\alpha} t^\alpha}\right) = e^{\frac{1}{\alpha} t^\alpha}.$

**Theorem 2.2.** Let $h(t) = f(g(t))$, $\alpha \in (0,1]$ and $f, g: (0, \infty) \to \mathbb{R}$ be $\alpha$-differentiable functions. Then $h(t)$ is $\alpha$-differentiable and for all $t$ with $t \neq 0$ and $g(t) \neq 0$ we have

$$T_\alpha(h)(t) = (T_\alpha f)(g(t))(T_\alpha g)(t) g(t)^{\alpha-1}. \qquad (3)$$

If $t = 0$, we have

$$T_\alpha(h)(t) = \lim_{t \to 0^+} (T_\alpha f)(g(t))(T_\alpha g)(t) g(t)^{\alpha-1}. \qquad (4)$$

**Theorem 2.3.** Let $f, g: [a, b] \to \mathbb{R}$ be two functions such that $f, g$ is differentiable. Then

$$\int_a^b f(t) T_\alpha(g)(t) d_\alpha(t) = fg|_a^b - \int_a^b g(t) T_\alpha(f)(t) d_\alpha(t). \qquad (5)$$

## 3. Conformable Fractional Integral

**Definition 3.1.** Let $f: [0, \infty) \to \mathbb{R}$. Then for all $t > 0$, $\alpha \in (0,1)$, the conformable fractional integral of $f$ of order $\alpha$ is

$$(I_\alpha f)(t) = \int_a^t f(x) d_\alpha(x) = \int_a^t x^{\alpha-1} f(x) dx \qquad (6)$$

where the integral is the usual Riemann improper integral.

**Theorem 3.1.** Let $f: [0, \infty) \to \mathbb{R}$ be any continuous function and $0 < \alpha \leq 1$. Then for all $t > a$,

$$T_\alpha I_\alpha f(t) = f(t). \qquad (7)$$

**Theorem 3.2.** Let $f: (a, b) \to R$ be differentiable and $0 < \alpha \leq 1$. Then for all $t > a$,

$$I_\alpha T_\alpha (f(t)) = f(t) - f(a) \qquad (8)$$

## 4. Existence and Uniqueness Theorems

Let linear sequential conformable fractional differential equation of order $n\alpha$

$$^n T_\alpha y + p_{n-1}(t)\, ^{n-1} T_\alpha y + \cdots + p_2(t)\, ^2 T_\alpha y + p_1(t) T_\alpha y + p_0(t) y = 0 \qquad (9)$$

where $^n T_\alpha y = T_\alpha T_\alpha \ldots T_\alpha y$. Analogously, non-homogeneous case is

$$^n T_\alpha y + p_{n-1}(t)\, ^{n-1} T_\alpha y + \cdots + p_2(t)\, ^2 T_\alpha y + p_1(t) T_\alpha y + p_0(t) y = f(t). \qquad (10)$$

We define an nth-order differential operator as following:

$$L_\alpha[y] = \,^n T_\alpha y + p_{n-1}(t)\, ^{n-1} T_\alpha y + \cdots + p_2(t)\, ^2 T_\alpha y + p_1(t) T_\alpha y + p_0(t) y = 0. \qquad (11)$$

**Theorem 4.1.**

Let $t^{\alpha-1}p(t), t^{\alpha-1}q(t) \in C(a,b)$ and let $y$ be $\alpha$-differentiable for $0 < \alpha \leq 1$. Then, the initial value problem

$$T_\alpha y + p(t)y = q(t) \tag{12}$$

$$y(t_0) = y_0 \tag{13}$$

has a unique solution on the interval $(a,b)$ where $t_0 \in (a,b)$.

**Proof**

We can write

$$T_\alpha y + p(t)y = q(t) \tag{14}$$

$$t^{1-\alpha}y' + p(t)y = q(t) \tag{15}$$

$$y' + t^{\alpha-1}p(t) = t^{\alpha-1}q(t). \tag{16}$$

From classical linear fundamental theorem existence and uniqueness, proof is clear.

**Theorem 4.2.**

Let $t^{\alpha-1}p_{n-1}(t), \ldots, t^{\alpha-1}p_1(t), t^{\alpha-1}p_0(t), t^{\alpha-1}q(t) \in C(a,b)$ and let $y$ be $n$ times $\alpha$-differentiable function. Then, a solution $y(t)$ of the initial-value problem

$$^nT_\alpha y + p_{n-1}(t)^{n-1}T_\alpha y + \cdots + p_2(t)\,^2T_\alpha y + p_1(t)T_\alpha y + p_0(t)y = q(t) \tag{17}$$

$$y(t_0) = y_0, T_\alpha y(t_0) = y_1, \ldots, {}^{n-1}T_\alpha y(t_0) = y_{n-1}, \ a < t_0 < b \tag{18}$$

exists on the interval and it is unique.

**Proof**

To show the existence of a local solution, we reduce our problem to the first order system of differential equations. So, let's make the following change of variables

$$x_1 = y, x_2 = T_\alpha y, x_3 = {}^2T_\alpha y, \ldots, x_n = {}^{n-1}T_\alpha y. \tag{19}$$

Hence, we have

$$\begin{aligned} T_\alpha x_1 &= x_2 \\ T_\alpha x_2 &= x_3 \\ &\vdots \\ T_\alpha x_{n-1} &= x_{n-1} \\ T_\alpha x_n &= -p_{n-1}x_n - \cdots - p_2 x_3 - p_1 x_2 - p_0 x_1 + q(t) \end{aligned} \tag{20}$$

Now, we can rewrite the problem in the following for

$$T_\alpha \begin{bmatrix} x_1 \\ x_2 \\ \vdots \\ x_{n-1} \\ x_n \end{bmatrix} + \underbrace{\begin{bmatrix} 0 & -1 & 0 & 0 & \cdots & 0 \\ 0 & 0 & -1 & 0 & \cdots & 0 \\ \vdots & \vdots & \vdots & \vdots & \cdots & \vdots \\ 0 & 0 & 0 & 0 & \cdots & -1 \\ p_0 & p_1 & p_2 & p_3 & \cdots & p_{n-1} \end{bmatrix}}_{P(t)} \underbrace{\begin{bmatrix} x_1 \\ x_2 \\ \vdots \\ x_{n-1} \\ x_n \end{bmatrix}}_{X(t)} = \underbrace{\begin{bmatrix} 0 \\ 0 \\ \vdots \\ 0 \\ q(t) \end{bmatrix}}_{Q(t)}, \quad (21)$$

$$T_\alpha X(t) + P(t)X(t) = Q(t), \quad (22)$$

$$X'(t) + t^{\alpha-1} P(t) X(t) = t^{\alpha-1} Q(t). \quad (23)$$

The existence and uniqueness of solution (17)-(18) follows from classical theorems on existence and uniqueness for system equation.

**Theorem 4.3.**

$L_\alpha$ is linear, i.e.

$$L_\alpha[c_1 y_1 + c_2 y_2] = c_1 L_\alpha[y_1] + c_2 L_\alpha[y_2] \quad (24)$$

where $y_1, y_2$ are $n$ times $\alpha$-differentiable functions and $c_1, c_2$ are arbitrary numbers.

**Proof**

From the definition conformable fractional derivative it follows that

$$\begin{aligned} L_\alpha[c_1 y_1 + c_2 y_2] &= {}^n T_\alpha (c_1 y_1 + c_2 y_2) + p_{n-1}(t)^{n-1} T_\alpha (c_1 y_1 + c_2 y_2) + \cdots \\ &\quad + p_1(t) T_\alpha (c_1 y_1 + c_2 y_2) + p_0(t)(c_1 y_1 + c_2 y_2) \\ &= c_1 \left( {}^n T_\alpha y_1 + p_{n-1}(t)^{n-1} T_\alpha y_1 + \cdots + p_0(t) y_1 \right) \\ &\quad + c_2 \left( {}^n T_\alpha y_2 + p_{n-1}(t)^{n-1} T_\alpha y_2 + \cdots + p_0(t) y_2 \right). \end{aligned}$$

$$(25)$$

So, we have

$$L_\alpha[c_1 y_1 + c_2 y_2] = c_1 L_\alpha[y_1] + c_2 L_\alpha[y_2].$$

**Theorem 4.4.**

Let $y_1, y_2, \ldots, y_n$ be the solutions of equation $L_\alpha[y] = 0$. In this case, the linear combination

$$y = c_1 y_1 + c_2 y_2 + \cdots + c_n y_n \quad (26)$$

is also its solution for the arbitrary constants $c_k, k = 1, \ldots, n$.

**Proof**

Let $y_1(t), y_2(t), \ldots, y_n(t)$ be solutions of the equation $L_\alpha y = 0$, for arbitrary constants $c_k, k = 1, \ldots, n$ and let us

$$y = c_1 y_1 + c_2 y_2 + \cdots + c_n y_n. \tag{27}$$

By using linear property of $L_\alpha$, we have

$$L_\alpha(y) = L_\alpha(c_1 y_1 + c_2 y_2 + \cdots + c_n y_n) = c_1 L_\alpha(y_1) + c_2 L_\alpha(y_2) + \cdots + c_n L_\alpha(y_n)$$
$$= 0 + \cdots + 0 = 0.$$

**Definition 4.1.**

We assume that $y_1(t), y_2(t), \ldots, y_n(t)$ are at least $(n-1)$ times $\alpha$-differentiable functions. For any $0 < \alpha \leq 1$, the determinant

$$W_\alpha(y_1, y_2, \ldots, y_n) = \begin{vmatrix} y_1 & y_2 & \cdots & y_n \\ T_\alpha y_1 & T_\alpha y_2 & \cdots & T_\alpha y_n \\ \vdots & \vdots & \cdots & \vdots \\ {}^{n-1}T_\alpha y_1 & {}^{n-1}T_\alpha y_2 & \cdots & {}^{n-1}T_\alpha y_n \end{vmatrix} \tag{28}$$

is called $\alpha$-Wronskian of the functions.

**Definition 4.2.**

Any set $y_1, y_2, \ldots, y_n$ of $n$ solutions of $L_\alpha[y] = 0$ is said to be a fundamental set of solutions.

**Theorem 4.5.**

Let $y_1, y_2, \ldots, y_n$ be $n$ solutions of $L_\alpha[y] = 0$. If there is a $t_0 \in (a, b)$ such that $W_\alpha(y_1, y_2, \ldots, y_n)(t_0) \neq 0$, then $\{y_1, y_2, \ldots, y_n\}$ is a fundamental set of solutions.

**Proof**

If $y(t)$ is a solution of $L_\alpha[y] = 0$, then we can write $y(t)$ as a linear combination of $y_1, y_2, \ldots, y_n$ which we need to show. That is

$$y = c_1 y_1 + c_2 y_2 + \cdots + c_n y_n. \tag{29}$$

So, the problem is now reduced to finding the constants $c_k, 1 \leq k \leq n$. We can write the following linear system of equations

$$\begin{array}{rcl}
c_1 y_1(t_0) + c_2 y_2(t_0) + \cdots + c_n y_n(t_0) & = & y(t_0), \\
c_1 T_\alpha y_1(t_0) + c_2 T_\alpha y_2(t_0) + \cdots + c_n T_\alpha y_n(t_0) & = & T_\alpha y(t_0), \\
\vdots & & \vdots \\
c_1 {}^{n-1}T_\alpha y_1(t_0) + c_2 {}^{n-1}T_\alpha y_2(t_0) + \cdots + c_n {}^{n-1}T_\alpha y_n(t_0) & = & {}^{n-1}T_\alpha y(t_0).
\end{array} \tag{30}$$

Applying Cramer's rule, we can get

$$c_k = \frac{W_\alpha^k(t_0)}{W_\alpha(t_0)}, \quad 1 \leq k \leq n. \tag{31}$$

Since $W_\alpha(t_0) \neq 0$, it follows that $c_1, c_2, \ldots, c_n$ exist.

**Theorem 4.6.**

Let $t^{\alpha-1}p_{n-1}(t),\ldots,t^{\alpha-1}p_1(t),t^{\alpha-1}p_0(t) \in C(a,b)$. Then, the equation $L_\alpha[y] = 0$ has a fundamental set of solutions $\{y_1, y_2, \ldots, y_n\}$.

**Proof**

Let $t_0 \in (a,b)$. Consider the following $n$ initial value problems

$$L_\alpha[y] = 0, y(t_0) = 1, T_\alpha y(t_0) = 0, \ldots, {}^{n-1}T_\alpha y(t_0) = 0$$
$$L_\alpha[y] = 0, y(t_0) = 0, T_\alpha y(t_0) = 1, \ldots, {}^{n-1}T_\alpha y(t_0) = 0$$
$$\vdots$$
$$L_\alpha[y] = 0, y(t_0) = 0, T_\alpha y(t_0) = 0, \ldots, {}^{n-1}T_\alpha y(t_0) = 1 \tag{32}$$

From Theorem 4.2, it follows that there is, for each index $i$, the solution $y_i$ of $ith$ problem. By the Theorem 4.5, the set of $\{y_1, y_2, \ldots, y_n\}$ is a fundamental set of solutions since

$$W_\alpha(t) = \begin{vmatrix} 1 & 0 & \cdots & 0 \\ 0 & 1 & \cdots & 0 \\ \vdots & \vdots & \cdots & \vdots \\ 0 & 0 & \cdots & 1 \end{vmatrix} = 1 \neq 0. \tag{33}$$

**Theorem 4.7.**

Let $y_1, y_2, \ldots, y_n$ be $n$ solutions of equation $L_\alpha[y] = 0$. Then

(1) $W_\alpha(t)$ satisfies the differential equation $T_\alpha W + p_{n-1}(t)W = 0$,
(2) For any $t_0 \in (a,b)$

$$W_\alpha(t) = W_\alpha(t_0) e^{-\int_{t_0}^{t} x^{\alpha-1}(p_{n-1}(x))dx} \tag{34}$$

Moreover, if $W_\alpha(t_0) \neq 0$ then $W_\alpha(t) \neq 0$ for all $t \in (a,b)$.

**Proof**

(1) Introducing new variables
$$x_1 = y, x_2 = T_\alpha y, x_3 = {}^2T_\alpha y, \ldots, x_n = {}^{n-1}T_\alpha y. \tag{35}$$

We can rewrite the differential in the matrix form

$$T_\alpha \underbrace{\begin{bmatrix} x_1 \\ x_2 \\ \vdots \\ x_{n-1} \\ x_n \end{bmatrix}}_{X(t)} = \underbrace{\begin{bmatrix} 0 & 1 & 0 & 0 & \cdots & 0 \\ 0 & 0 & 1 & 0 & \cdots & 0 \\ \vdots & \vdots & \vdots & \vdots & \cdots & \vdots \\ 0 & 0 & 0 & 0 & \cdots & 1 \\ -p_0 & -p_1 & -p_2 & -p_3 & \cdots & -p_{n-1} \end{bmatrix}}_{P(t)} \underbrace{\begin{bmatrix} x_1 \\ x_2 \\ \vdots \\ x_{n-1} \\ x_n \end{bmatrix}}_{X(t)} \tag{36}$$

$$T_\alpha X(t) = P(t)X(t) \tag{37}$$

Therefore, we obtain

$$T_\alpha W_\alpha(t) = (P_{11} + P_{22} + \cdots + P_{nn})W_\alpha(t). \tag{38}$$

Really, if α- conformable derivative of $W_\alpha(y_1, y_2, \ldots, y_n)$ is calculated, then we get

$$T_\alpha W_\alpha(t) = -p_{n-1}(t)W_\alpha(t). \tag{39}$$

Hence

$$\frac{T_\alpha(W_\alpha(t))}{W_\alpha(t)} = -p_{n-1}(t) \tag{40}$$

$$\ln(W_\alpha(t)) - \ln(W_\alpha(t_0)) = -\int_{t_0}^{t} x^{\alpha-1} p_{n-1}(x) dx \tag{41}$$

$$W_\alpha(t) = W_\alpha(t_0) e^{-\int_{t_0}^{t} x^{\alpha-1} p_{n-1}(x) dx}. \tag{42}$$

So, proof of theorem is finished.

**Theorem 4.8.**

Let $t^{\alpha-1} p_{n-1}(t), \ldots, t^{\alpha-1} p_1(t), t^{\alpha-1} p_0(t) \in C(a,b)$. If $\{y_1, y_2, \ldots, y_n\}$ is a fundamental set of solutions of $L_\alpha[y] = 0$, then $W_\alpha(t) \neq 0$ for all $t \in (a,b)$.

**Proof**

Suppose that $t_0$ be any point in $(a,b)$. By Theorem 4.2, there is a unique solution $y(t)$ of the initial value problem

$$L_\alpha[y] = 0, y(t_0) = 1, T_\alpha y(t_0) = 0, \ldots, {}^{n-1}T_\alpha y(t_0) = 0. \tag{43}$$

There exist unique constants $c_1, c_2, \ldots, c_n$ such that

$$\begin{array}{rcl}
c_1 y_1(t) + c_2 y_2(t) + \cdots + c_n y_n(t) & = & y(t) \\
c_1 T_\alpha y_1(t) + c_2 T_\alpha y_2(t) + \cdots + c_n T_\alpha y_n(t) & = & T_\alpha y(t) \\
\vdots & & \vdots \\
c_1 {}^{n-1}T_\alpha y_1(t) + c_2 {}^{n-1}T_\alpha y_2(t) + \cdots + c_n {}^{n-1}T_\alpha y_n(t) & = & {}^{n-1}T_\alpha y(t)
\end{array} \tag{44}$$

for all $t \in (a,b)$ since $\{y_1, y_2, \ldots, y_n\}$ is a fundamental set of solutions. In particular, for $t = t_0$, we obtain the system

$$\begin{array}{rcl}
c_1 y_1(t_0) + c_2 y_2(t_0) + \cdots + c_n y_n(t_0) & = & 1 \\
c_1 T_\alpha y_1(t_0) + c_2 T_\alpha y_2(t_0) + \cdots + c_n T_\alpha y_n(t_0) & = & 0 \\
\vdots & & \vdots \\
c_1 {}^{n-1}T_\alpha y_1(t_0) + c_2 {}^{n-1}T_\alpha y_2(t_0) + \cdots + c_n {}^{n-1}T_\alpha y_n(t_0) & = & 0
\end{array} \tag{45}$$

This system has a unique solution

$$c_k = \frac{W_\alpha^k}{W_\alpha(t_0)}, \quad 1 \leq k \leq n. \tag{46}$$

Here, for every $k$

$$W_\alpha^k = \begin{bmatrix} y_1(t_0) & y_2(t_0) & \cdots & y_{k-1}(t_0) & 1 & y_k(t_0) & \cdots & y_n(t_0) \\ T_\alpha y_1(t_0) & T_\alpha y_2(t_0) & \cdots & T_\alpha y_{k-1}(t_0) & 0 & T_\alpha y_k(t_0) & \cdots & T_\alpha y_n(t_0) \\ \vdots & \vdots & \vdots & \vdots & \vdots & \vdots & \vdots & \vdots \\ {}^{n-1}T_\alpha y_1(t_0) & {}^{n-1}T_\alpha y_2(t_0) & \cdots & {}^{n-1}T_\alpha y_{k-1}(t_0) & 0 & {}^{n-1}T_\alpha y_k(t_0) & \cdots & {}^{n-1}T_\alpha y_n(t_0) \end{bmatrix}$$

(47)

For existence $c_1, c_2, \ldots, c_n$, it should be $W_\alpha(t_0) \neq 0$. By Theorem 4.7, we conclude that $W_\alpha(t_0) \neq 0$ for all $t_0 \in (a, b)$.

**Theorem 4.9.**

Let $t^{\alpha-1}p_{n-1}(t), \ldots, t^{\alpha-1}p_1(t), t^{\alpha-1}p_0(t) \in C(a, b)$. The solution set $\{y_1, y_2, \ldots, y_n\}$ is a fundamental set of solutions of equation $L_\alpha[y] = 0$ if and only if the functions $y_1, y_2, \ldots, y_n$ are linearly independent.

**Proof**

Suppose that $\{y_1, y_2, \ldots, y_n\}$ is a fundamental set of solutions, then by Theorem 4.8 it follows that there is $t_0 \in (a, b)$ such that $W_\alpha(t_0) \neq 0$. Assume that

$$c_1 y_1(t) + c_2 y_2(t) + \cdots + c_n y_n(t) = 0 \tag{48}$$

for all $t \in (a, b)$. By repeated $\alpha$-differentiation of the equation (48), we find

$$\begin{array}{rcl} c_1 y_1(t) + c_2 y_2(t) + \cdots + c_n y_n(t) &=& 0 \\ c_1 T_\alpha y_1(t) + c_2 T_\alpha y_2(t) + \cdots + c_n T_\alpha y_n(t) &=& 0 \\ \vdots & & \vdots \\ c_1{}^{n-1}T_\alpha y_1(t) + c_2{}^{n-1}T_\alpha y_2(t) + \cdots + c_n{}^{n-1}T_\alpha y_n(t) &=& 0 \end{array} \tag{49}$$

Thus, one finds $c_1, c_2, \ldots, c_n$ by solving the system

$$\begin{array}{rcl} c_1 y_1(t_0) + c_2 y_2(t_0) + \cdots + c_n y_n(t_0) &=& 0 \\ c_1 T_\alpha y_1(t_0) + c_2 T_\alpha y_2(t_0) + \cdots + c_n T_\alpha y_n(t_0) &=& 0 \\ \vdots & & \vdots \\ c_1{}^{n-1}T_\alpha y_1(t_0) + c_2{}^{n-1}T_\alpha y_2(t_0) + \cdots + c_n{}^{n-1}T_\alpha y_n(t_0) &=& 0 \end{array} \tag{50}$$

Namely, by using Cramer's rule one finds

$$c_1 = c_2 = \cdots = c_n = \frac{0}{W_\alpha(t_0)} = 0. \tag{51}$$

Thus, set of functions $\{y_1, y_2, \ldots, y_n\}$ are linearly independent.

Contrarily, suppose $\{y_1, y_2, \ldots, y_n\}$ is a linearly independent set and suppose $\{y_1, y_2, \ldots, y_n\}$ is not a fundamental set of solutions. Then, by Theorem 4.5, we get $W_\alpha(t) = 0$, for all $t \in (a, b)$. We choose any $t_0 \in (a, b)$. Then $W_\alpha(t_0) = 0$. But this says $W_\alpha(t_0) \neq 0$ that the matrix

$$\begin{bmatrix} y_1(t_0) & y_2(t_0) & \cdots & y_n(t_0) \\ T_\alpha y_1(t_0) & T_\alpha y_2(t_0) & \cdots & T_\alpha y_n(t_0) \\ \vdots & \vdots & \cdots & \vdots \\ {}^{n-1}T_\alpha y_1(t_0) & {}^{n-1}T_\alpha y_2(t_0) & \cdots & {}^{n-1}T_\alpha y_n(t_0) \end{bmatrix} \tag{52}$$

is not invertible which means that there exist $c_1, c_2, \ldots, c_n$, $c_1^2 + c_2^2 + \cdots + c_n^2 \neq 0$ such that

$$\begin{aligned} c_1 y_1(t_0) + c_2 y_2(t_0) + \cdots + c_n y_n(t_0) &= 0 \\ c_1 T_\alpha y_1(t_0) + c_2 T_\alpha y_2(t_0) + \cdots + c_n T_\alpha y_n(t_0) &= 0 \\ &\vdots \quad \vdots \\ c_1{}^{n-1}T_\alpha y_1(t_0) + c_2{}^{n-1}T_\alpha y_2(t_0) + \cdots + c_n{}^{n-1}T_\alpha y_n(t_0) &= 0 \end{aligned} \tag{53}$$

Now, let

$$y(t) = c_1 y_1(t) + c_2 y_2(t) + \cdots + c_n y_n(t) \tag{54}$$

for all $t \in (a, b)$. Then, $y(t)$ is the solution of the differential equation and

$$y(t_0) = T_\alpha y(t_0) = \cdots = {}^{n-1}T_\alpha y(t_0) = 0. \tag{55}$$

But the zero function also is the solution of the initial value problem. By Theorem 4.2, it should be

$$c_1 y_1(t) + c_2 y_2(t) + \cdots + c_n y_n(t) = 0 \tag{56}$$

for all $t \in (a, b)$ with $c_1, c_2, \ldots, c_n$ not all equal to zero which means that $y_1, y_2, \ldots, y_n$ are linearly dependent which is not in agreement with our earlier assumption that $y_1, y_2, \ldots, y_n$ are linearly independent.

**Theorem 4.10.**

Let $y_1, y_2, \ldots, y_n$ be linearly independent solutions of the equation $L_\alpha[y] = 0$. Then the general solution of the equation is

$$y = c_1 y_1 + c_2 y_2 + \cdots + c_n y_n \tag{57}$$

where the arbitrary constants $c_k, k = 1, \ldots, n$.

**Proof**

Particular solution at any $t = t_0$ is obtained by the help of initial conditions as following:

$$y(t_0) = \gamma_0, T_\alpha y(t_0) = \gamma_1, \ldots, {}^{n-1}T_\alpha y(t_0) = \gamma_{n-1} \tag{58}$$

where $t_0 \in (a, b)$ and $y_0, y_1, \ldots, y_{n-1}$ are arbitrary constants. If we choose $c_1, c_2, \ldots, c_n$ constants to protect conditions (58), then proof is completed. To accomplish this, we can write following system of equations:

$$c_1 y_1(t_0) + c_2 y_2(t_0) + \cdots + c_n y_n(t_0) = \gamma_0$$
$$c_1 T_\alpha y_1(t_0) + c_2 T_\alpha y_2(t_0) + \cdots + c_n T_\alpha y_n(t_0) = \gamma_1$$
$$\vdots \qquad (59)$$
$$c_1{}^{n-1}T_\alpha y_1(t_0) + c_2{}^{n-1}T_\alpha y_2(t_0) + \cdots + c_n{}^{n-1}T_\alpha y_n(t_0) = \gamma_{n-1}$$

Writing the above system in matrix form, we get

$$\begin{bmatrix} y_1(t_0) & y_2(t_0) & \cdots & y_n(t_0) \\ T_\alpha y_1(t_0) & T_\alpha y_2(t_0) & \cdots & T_\alpha y_n(t_0) \\ \vdots & \vdots & \vdots & \vdots \\ {}^{n-1}T_\alpha y_1(t_0) & {}^{n-1}T_\alpha y_2(t_0) & \cdots & {}^{n-1}T_\alpha y_n(t_0) \end{bmatrix} \begin{bmatrix} c_1 \\ c_2 \\ \vdots \\ c_n \end{bmatrix} = \begin{bmatrix} \gamma_0 \\ \gamma_1 \\ \vdots \\ \gamma_{n-1} \end{bmatrix} \qquad (60)$$

Because $y_1, y_2, \ldots, y_n$ is linearly independent solutions of $L_\alpha[y] = 0$, $W_\alpha(t_0) \neq 0$. In this case, according to the fundamental theorem of algebra, system (59) has a unique solution. Hence, proof is completed.

**Theorem 4.11.**

Let $y_p$ be any particular solution of the nonhomogeneous linear nth-order sequential differential equation (10) and $\{y_1, y_2, \ldots, y_n\}$ be a fundamental set of solutions of the associated homogeneous differential equation (9). Then the general solution of the equation is

$$y = c_1 y_1 + c_2 y_2 + \cdots + c_n y_n + y_p \qquad (61)$$

where the arbitrary constants $c_k, k = 1, \ldots, n$.

**Proof**

Let $L_\alpha$ be the differential operator and let $Y(t)$ and $y_p(t)$ be particular solutions of the nonhomogeneous equation $L_\alpha[y] = q(t)$. If we define $u(t) = Y(t) - y_p(t)$, then by linearity of $L_\alpha$ we have

$$L_\alpha[u] = L_\alpha[Y(t) - y_p(t)] = L_\alpha[Y(t)] - L_\alpha[y_p(t)] = q(t) - q(t) = 0. \qquad (62)$$

This shows that $u(t)$ is a solution of the homogeneous equation $L_\alpha[y] = 0$. Hence, by Theorem 4.4,

$$u(t) = c_1 y_1(t) + c_2 y_2(t) + \cdots + c_n y_n(t) \qquad (63)$$

and so

$$Y(t) - y_p(t) = c_1 y_1(t) + c_2 y_2(t) + \cdots + c_n y_n(t) \qquad (64)$$

or

$$Y(t) = c_1 y_1(t) + c_2 y_2(t) + \cdots + c_n y_n(t) + y_p(t). \qquad (65)$$

## 5. Conclusion

In our work, we present existence and uniqueness theorems for sequential linear conformable fractional differential equations. It has been found that results obtained from this work is analogous to the results obtained from the ordinary case.